\documentclass[draft]{amsart}
\usepackage{times}
\usepackage{amsfonts}
\usepackage{amsmath}
\usepackage{amscd}
\usepackage{amssymb}
\usepackage{latexsym}
\usepackage{natbib}

\theoremstyle{plain}

\swapnumbers
\newtheorem{theorem}[subsection]{Theorem}
\newtheorem{corollary}[subsection]{Corollary}
\newtheorem{lemma}[subsection]{Lemma}
\newtheorem{proposition}[subsection]{Proposition}

\newtheorem*{maintheorem}{Main Theorem}

\newtheorem*{citedtheorem}{Theorem}

\theoremstyle{definition}
\newtheorem{definition}[subsection]{Definition}

\newtheorem{example}[subsection]{Example}

\theoremstyle{remark}

\newtheorem{remark}[subsection]{Remark}

\swapnumbers

\newcommand\myhomepage{http://www.math.ohio-state.edu/\-\~{}schoutens}

\newcommand{\emptyprop}{q}

\newcommand \zet{\mathbb Z}
\newcommand \numgen[2]{\textsf{\textbf{#1}}(#2)}
\newcommand \nd{non-degenerate}

\newcommand \complet[1]{\widehat {#1}} 
\newcommand \exactseq [5]{0\to{#1}\:\map{#2}\:{#3}\:\map{#4}\:{#5}\to0}
\newcommand \Exactseq [3]{0\to {#1}\to {#2}\to {#3}\to 0}	
\newcommand \ext[4]{\operatorname{Ext}_{#1}^{#2}(#3,#4)}
 
\newcommand \id{\mathfrak a}

\newcommand \inverse[2]{{#1^{-1}(#2)}}
\newcommand \iso{\cong}

\newcommand \map[1]{{\newcommand{\tmpprop}{#1q}  \if\tmpprop\emptyprop \to\else \xrightarrow{{\phantom{i}{#1}\phantom{i}}}\fi}} 
\newcommand \maxim{\mathfrak m}

\newcommand \onto{\twoheadrightarrow}
\newcommand \pol[2]{#1[#2]}
\newcommand \pow[2]{#1[[#2]]}
\newcommand \pr{\mathfrak p}
\newcommand \PR{\mathfrak P}

\newcommand \rij[2]{(#1_1,\dots,#1_{#2})}

\newcommand \tensor{\otimes}

\newcommand \op\operatorname

\newcommand \bound[2]{\textsf{#1}(#2)}

\newcommand \UFD{unique factorization domain}
\newcommand \ch{characteristic}
\newcommand \homo{homomorphism}
\newcommand \CM{Coh\-en-Mac\-au\-lay}
\renewcommand\iff{if, and only if,}

\title[Number of generators of a Cohen-Macaulay ideal] {Absolute bounds on the number of generators of Cohen-Macaulay ideals of height at most $2$}
\author{Hans Schoutens
}
\thanks{Partially supported by a  grant from the National Science Foundation.}\date{03.04.2003}
\address{Department of Mathematics\\
Ohio State University\\
Columbus, OH 43210 (USA)}
\email{schoutens@math.ohio-state.edu}
\urladdr{\myhomepage} 
\subjclass{13E15,13C14}

\begin{document}

\begin{abstract}
For a Noetherian local domain $A$, there exists an upper bound $N_\tau(A)$ on the minimal number of generators of any height two ideal $\id$ for which $A/\id$ is \CM\ of type $\tau$.  More precisely, we may take $N_\tau(A):=(\tau+1)e_{\text{h}}(A)$, where $e_{\text{h}}(A)$ is the homological multiplicity of $A$.
\end{abstract}

\keywords{number of generators, Cohen-Macaulay ideals, Noether Normalization, homological multiplicity}

\maketitle


\section{Introduction}

In this paper, we study the minimal number of generators of certain classes of ideals. Namely, let $\mathcal I(A)$ be a class of ideals in a Noetherian ring $A$. Let $\numgen iA$ be the maximum  of all $\mu_A(\id)$ with $\id\in\mathcal I(A)$ (we allow $\numgen iA=\infty$), where in general,  $\mu_A(M)$ denotes the minimal number of generators of a finitely generated $A$-module $M$. Let us call $\numgen iA$ the \emph{girth} of the class $\mathcal I(A)$. One can pose the following two questions  about these girths.
\begin{itemize}
\item For which classes of ideals $\mathcal I(A)$ is its girth $\numgen iA$ finite? 
\item If the girth $\numgen iA$ is finite, can it be described, or at least be bounded, in terms of other invariants of $A$? 
\end{itemize}

For instance, \citet{SV}  use Noether Normalization to show that any prime ideal in a  two-dimensional affine algebra $A$ (that is to say, a two-dimensional finitely generated algebra over a field) is generated by at most $N(A)$ elements, where $N(A)$ only depends on the algebra. In our terminology, the girth of the prime spectrum of a two-dimensional affine algebra is bounded. In this paper, we will address both questions for the girths of the following classes of ideals: (1) the girth $\textsf{cm}_1(A)$ of the class of height one \CM\ ideals (that is to say, the height one ideals $\id$ for which $A/\id$ is \CM), (2) the girth $\textsf{g}_2(A)$ of the class of height two Gorenstein ideals (that is to say, the height two ideals $\id$ for which $A/\id$ is Gorenstein), and, more generally, (3) the girth $\textsf{cm}_2^\tau(A)$ of the class of height two \CM\ ideals for which $A/\id$ has type $\tau$.  Note that $\textsf{g}_2(A)=\textsf{cm}_2^1(A)$. To state precise results, we need a definition.

Let $A$ be a Noetherian local ring with residue field $k$. We call $A$ \emph{\nd} if $A$ has the same \ch\ as any of its irreducible components, that is to say, $\op{char}(A)=\op{char}(A/\pr)$, for every minimal prime $\pr$ of $A$. Note that this condition is void if $A$ is equi\ch. In mixed \ch, it means that $A$ has \ch\ zero and its residue field has \ch\ $p$, and, moreover,  $p$ does not belong to any minimal prime of $A$. In particular, any Noetherian local domain is \nd. Observe that the class of \nd\ local rings is closed under completion: this is clear for equi\ch\ rings; for mixed \ch\ rings, if $p$ lies in a minimal prime $\pr$ of the completion of $A$, then it lies in $\pr\cap A$, which is necessarily a minimal prime of $A$. The motivation for introducing this class of local rings comes from the following structure theorem due to Cohen (see \citep[Theorem 29.4 and Remark]{Mats} or \citep[IX. Th\'eor\`eme 3]{Bour83}).

\begin{citedtheorem}[Cohen Structure Theorem]
Let $A$ be a complete Noetherian local ring. If $A$ is \nd, then it admits a regular (complete local) subring $S$ over which it is finite.
\end{citedtheorem}

Combining Theorems~\ref{T:locCM} and \ref{T:locGor} and Corollaries~\ref{C:locCM2} and \ref{C:hommult} gives the main result of this paper. 

\begin{maintheorem}
Let $A$ be a \nd\ Noetherian local ring. Each of the following are finite $\textsf{cm}_1(A)$,  $ \textsf{g}_2(A)$ and, more generally,  $\textsf{cm}^\tau_2(A)$ (for fixed $\tau$). 

In fact, these girths are bounded respectively by $e$, $2e$ and $(\tau+1)e$ where, in general, $e$ is the homological multiplicity of $A$, but in case $A$ is equi\ch\ and \CM, we may take $e$ to be the (usual) multiplicity of $A$.
\end{maintheorem}

The bounds in case $A$ is \CM\ have been obtained already by Rees   and Sally respectively, without any equi\ch\ assumption; see \citep[Chapter V, Theorem 3.2 and Corollary 3.3]{Sally}. Bounds on the minimal number of generators of a \CM\ ideal $\id$ of arbitrary height or without the \CM\ assumption, but involving invariants of the residue ring $A/\id$ (whence not providing yet the finiteness of the above girths), can be found in   \citep{BER,Sally,SV,SchCMId,VallaGen,VasCD, VasMultGen}. For the height two case, the bounds in the Main Theorem are an improvement of the ones given in \citep[Example 9.5.1]{Vas}.

Note that the example $A:=\pow K{X,Y}/(X^2,XY)$ and the ideal $\id:=(X,Y)A$ shows that when $A$ is not \CM, then we cannot take in the Main Theorem $e$ equal to the multiplicity of $A$; in fact, in this example the stated bound is optimal: $A$ has multiplicity $1$ and homological multiplicity $2$. 

The proof of the Main Theorem uses a technique due to  \citet{SV}, except that we replace their use of Noether Normalization by the Cohen Structure Theorem. Using the Forster-Swan Theorem, we obtain estimates $\textsf{cm}_1(A)\leq e+d$ and $\textsf{cm}^\tau_2(A)\leq (\tau+1)e+d$, for any $d$-dimensional   Noetherian  ring $A$ which is finitely generated by at most $e$ elements over some regular subring. A particular instance of these global bounds is given by the following theorem.

\begin{theorem}\label{T:A}
Let $Y\to X$ be a finite dominant map of degree $f$ between affine $d$-dimensional schemes. If $X$ has no singularities, then every Gorenstein subscheme $W$ of $Y$ of codimension $2$ is the (ideal-theoretic) intersection of at most $2f+3d-2$ hypersurfaces. 

If $Y$ is moreover \CM, we require at most $2f+d-2$ hypersurfaces.
\end{theorem}

In the last section, bounds for  affine algebras are shown to be uniform, in the sense that the bounds on the girths only depend on the degree of the polynomials representing the affine algebra as a homomorphic image of a polynomial ring (see Theorems~\ref{T:eqch} and \ref{T:aff}). 

\subsubsection*{Acknowledgement}
This paper arose in part through several useful and stimulating conversations I had with  W. Vasconcelos.

\section{The Height One Case}

\begin{proposition}\label{P:finCM}
Let $A\to B$ be a local and finite \homo\ of Noetherian local rings. If $B$ is a \CM\ local ring, then it is also \CM\ as an $A$-module.
\end{proposition}
\begin{proof}
We will induct on the Krull dimension $d$ of $B$. There is nothing to prove if $d=0$, so assume $d>0$. I claim that $B$ has positive depth as an $A$-module. Indeed, if not, then the maximal ideal $\maxim$ of $A$ is an associated prime of $B$ viewed as an $A$-module. In other words, there is some non-zero $s\in B$, such that $s\maxim=0$ in $B$. Since $B$ is finite over $A$, the ideal $\maxim B$ is $\mathfrak n$-primary, where $\mathfrak n$ is the maximal ideal of $B$. Therefore, some power $\mathfrak n^e$ is contained in $\maxim B$. It follows that $s\mathfrak n^e=0$, showing that $\mathfrak n$ is an associated prime of $B$. However, this contradicts our assumption that $B$ has depth $d>0$, thus establishing  the claim.  Therefore, we can find $r\in \maxim$, such that $r$ is $B$-regular. It follows that $B/rB$ is a local \CM\ ring of dimension $d-1$. By induction, $B/rB$ is \CM\ when viewed as an $A/rA$-module.  It follows that $B$ is \CM\ as an $A$-module, as required.
\end{proof}

To reduce to the complete case, we need the following easy lemma.

\begin{lemma}\label{L:ff}
Let $A\to B$ be a faithfully flat \homo\ between Noetherian local rings. For an arbitrary ideal of $A$, both its minimal number of generators and its height  remain the same after extension to $B$.
\end{lemma}
\begin{proof}
Let $\id$ be an ideal of $A$ and $h$ its height. Let $\pr$ be a minimal prime of $\id$ of height $h$. By faithful flatness, the fiber ring $B_\pr/\pr B_\pr$ is non-zero, whence admits a minimal prime $\PR$.  In particular, $\pr=\PR\cap A$. By \citep[Theorem 15.1]{Mats}, flatness then yields that   $\op{ht}(\pr)=\op{ht}(\PR)$. Since $\id B$ is contained in $\PR$,  its height is at most $h$. On the other hand, if its height were to  be less than $h$, then there would exist a minimal prime ideal $\mathfrak Q$ of $\id B$ of height less than $h$, whose contraction to $A$ would be a prime ideal of height less than $h$ (again using \citep[Theorem 15.1]{Mats}) and containing $\id$, contradiction.

As for the minimal number of generators, let $n:=\mu_A(\id)$ and choose minimal generators $x_1,\dots,x_n$ for $\id$. By  the Nakayama  Lemma we can renumber the $x_i$ so that $x_1,\dots, x_m$ generate $\id B$ minimally.  By faithfully flatness, we get 
	\begin{equation*}
	\id=\id B\cap A= \rij xmB\cap A=\rij xmA,
	\end{equation*}
showing that $m=n$.
\end{proof}

\begin{theorem}\label{T:locCM}
For a \nd\  Noetherian local ring $A$, the girth invariant $\textsf{cm}_1(A)$ is finite, that is to say, there is an upper bound on the number of generators of an arbitrary height one \CM\ ideal of $A$.
\end{theorem}
\begin{proof}
By Lemma~\ref{L:ff}, we may assume that $A$ is moreover complete (note that the completion of a \nd\ ring is again \nd). By the Cohen Structure Theorem, there exists a regular local subring $S$ of $A$, such that $A$ is module finite over $S$. In particular, there exists a surjective \homo
	\begin{equation}
	\pi\colon S^N\onto A.
	\end{equation}
We will show that $\textsf{cm}_1(A)\leq N$. To this end, let $\id$ be an arbitrary height one \CM\ ideal in $A$. Let $d$ be the dimension of $A$, which is then also the dimension of $S$. From now on, we will view $A/\id$ only as an $S$-module. As such, its dimension is equal to the dimension of $S/(\id\cap S)$. Since $S\subset A$ is finite, $\id\cap S$ has height one, whence $S/\id\cap S$ has dimension $d-1$.  Together with Proposition~\ref{P:finCM}, we get that $A/\id$ is a \CM\ $S$-module of dimension $d-1$. In particular, its depth is $d-1$, whence, by the Auslander-Buchsbaum Theorem, its projective dimension is $1$. Thus from the exact sequence
	\begin{equation}
	\Exactseq {\inverse\pi\id} {S^N} {A/\id}
	\end{equation}
we get that $\inverse\pi\id$ is a free $S$-module, necessarily therefore of rank at most $N$. Since $\inverse\pi\id$ is generated by at most $N$ elements, so is $\id$.
\end{proof}

\section{The  Height Two Case}

\begin{theorem}\label{T:locGor}
For a \nd\ Noetherian local ring $A$, the girth invariant $\textsf{g}_2(A)$ is finite, that is to say, there is an upper bound on the number of generators of an arbitrary height two  ideal $\id$ of $A$ for which $A/\id$ is Gorenstein.
\end{theorem}
\begin{proof}
As in the proof of Theorem~\ref{T:locCM}, we may assume that $A$ is moreover complete and finitely generated as a module over a $d$-dimensional regular local subring $S$. Choose some  surjective \homo\ $\varphi\colon S^N\to A$ and let $W:=\inverse\varphi\id$. In particular, we have an exact sequence
	\begin{equation}\label{eq:w}
	\Exactseq W{S^N}{A/\id}.
	\end{equation} 
As before, since this time $A/\id$ is  \CM\ $S$-module of dimension $d-2$, it follows from Proposition~\ref{P:finCM} and the Auslander-Buchsbaum Theorem, that $A/\id$ has projective dimension $2$ as an $S$-module. Therefore, $W$ has projective dimension $1$. Let
	\begin{equation}\label{eq:res}
	\exactseq{S^p}f{S^q}{}W
	\end{equation}
be a minimal free $S$-resolution of $W$, so that $W$ is minimally generated by $q$ elements. Taking the ($S$-)dual of sequence~\eqref{eq:res}, we get an exact sequence
	\begin{equation}\label{eq:dual}
	S^q\map{f^*} S^p \to \ext S1WS\to 0,
	\end{equation}
where $f^*$ is the \emph{transpose} of $f$, that is to say, if $\mathbb A$ is a matrix defining $f$, then the matrix transpose of $\mathbb A$ gives the \homo\ $f^*$. In particular, since we took \eqref{eq:res} to be minimal, $\mathbb A$ has all its entries in the maximal ideal of $S$. Therefore, the same is true for $f^*$, so that by the Nakayama  Lemma, $\ext S1WS$ is minimally generated by $p$ elements.

Applying \citep[Theorem 3.3.7.(b)]{BH} to the finite local \homo\ $S\to A/\id$, we get that 
	\begin{equation*}
	\ext S2{A/\id}S=\omega_{A/\id},
	\end{equation*} 
where $\omega_{A/\id}$ is the canonical module of $A/\id$. However, since $A/\id$ is Gorenstein, we have that $\omega_{A/\id}\iso A/\id$. On the other hand, taking the dual of the exact sequence~\eqref{eq:w} shows that $\ext S1WS\iso \ext S2{A/\id}S$. In summary, we obtain an isomorphism
	\begin{equation*}
	\ext S1WS\iso A/\id.
	\end{equation*} 
Since this $S$-module is minimally generated by $p$ elements, \eqref{eq:w} yields that $p\leq N$. From the exact sequence
	\begin{equation*}
	0\to S^p\to S^q\to S^N
	\end{equation*}
we get that $q\leq p+N$. Therefore, $q\leq 2N$, showing that $W$, and hence a fortiori $\id$, can be generated by at most $2N$ elements.
\end{proof}

\begin{corollary}\label{C:locCM2}
For a \nd\ Noetherian local ring $A$ and an arbitrary $\tau\geq 1$, the girth invariant $\textsf{cm}_2^\tau(A)$ is finite.
\end{corollary}
\begin{proof}
Analyzing the proof of Theorem~\ref{T:locGor}, we see that the only place were we used that $A/\id$ is Gorenstein, is to establish the isomorphism $\omega_{A/\id}\iso A/\id$. If $A/\id$ is merely \CM\ of type $\tau$, then the canonical module $\omega_{A/\id}$ is generated as an $A/\id$-module by $\tau$ elements (\citep[Proposition 3.3.11]{BH}). Therefore, there is an epimorphism $(A/\id)^\tau\onto \omega_{A/\id}$. If $A$ is generated as an $S$-module by $N$ elements, then this implies that $\omega_{A/\id}$ is generated by at most $\tau N$ elements as an $S$-module. Hence from  \eqref{eq:w} and \eqref{eq:dual} we get that  $p\leq \tau N$ (notation as in that proof), so that $\mu_A(\id)\leq (\tau+1)N$.
\end{proof}

\section{Noether Normalization Degree}

From the above proofs it has become clear that the following invariant can be used to bound the minimal number of generators of an ideal.

\begin{definition}[Noether Normalization Degree]
Let $A$ be a Noetherian ring. We call its \emph{Noether Normalization degree} the least possible value of $\mu_S(A)$, where $S$ runs over all regular subrings of $A$ (this includes the case that there is no such regular subring over which $A$ is finite, in which case we set its Noether Normalization degree equal to $\infty$).
\end{definition}

A regular subring $S$ over which $A$ is finite is sometimes called a \emph{Cohen} ring. By the classical Noether Normalization Theorem, any finitely generated algebra over a field has finite Noether Normalization degree. By the Cohen Structure Theorem, so does any \nd\ complete Noetherian local ring.

\begin{corollary}\label{C:NNdeg}
Let $A$ be a Noetherian local ring and $N$ the Noether Normalization degree of its completion. We have estimates $\textsf{cm}_1(A)\leq N$ and $\textsf{cm}^\tau_2(A)\leq (\tau+1)N$.
\end{corollary}

Clearly, a Noetherian ring $A$ is regular \iff\ its Noether Normalization degree is equal to one.  More generally, we have the following result linking Noether Normalization degree and multiplicity. 

\begin{proposition}\label{P:NNmult}
For $A$ an equi\ch\ complete \CM\ ring with infinite residue field, its Noether Normalization degree is equal to its multiplicity.
\end{proposition}
\begin{proof}
Suppose $A$ has dimension $d$. Let $\maxim$ be its maximal ideal and $k$ its residue field. Let $N$ and $e$ be the respective Noether Normalization degree and the multiplicity of $A$. By the Cohen Structure Theorem, $k$ can be viewed as a subfield of $A$. By \citep[Theorem 14.14]{Mats}, we can choose $d$ elements $x_i\in\maxim$ which generate a reduction of $\maxim$. By \citep[Theorem 17.11]{Mats},  the $k$-vector space $A/\rij xdA$ has dimension $e$. The argument in \citep[Theorem 29.4.(iii)]{Mats} shows that $S:=\pow k{x_1,\dots,x_d}$ is a power series ring whence in particular a regular local ring. Moreover, by \citep[Theorem 8.4]{Mats}, if we choose liftings $y_1,\dots,y_e$ in $A$ of a basis of $A/\rij xdA$ over $k$, then the $y_i$ generate $A$ as an $S$-module. This shows that $N\leq e$. 

Conversely, assume $S$ is a regular subring of $A$ such that $\mu_S(A)=N$. Let $\rij xd$ be a regular sequence of parameters in $S$. It follows that $\rij xd$ is a system of parameters of $A$. By the Nakayama  Lemma, $A/\rij xdA$ has dimension $N$ over $k$, so that the ideal $\rij xdA$ has multiplicity $N$ by another application of \citep[Theorem 17.11]{Mats}. Since $e$ is the multiplicity of $\maxim$, we get from \citep[Formula 14.4]{Mats} that $e\leq N$.
\end{proof}

Note, for the last inequality, $e\leq N$, we did not use the \CM\ nor the equi\ch\ assumption (use \citep[Theorem 14.10]{Mats}), so that the Noether Normalization degree is always at least as big as the multiplicity of $A$. The example $A:=\pow k{X,Y}/(X^2,XY)$ shows that this inequality is in general strict (here $e=1$ but $N=2$). In fact, $\pow{\zet_p}X/(pX)$, where $\zet_p$ is the ring of $p$-adic integers, is an example of a complete Noetherian local ring with infinite Noether Normalization degree but multiplicity $e=2$. From the proof, it also follows that there is a close connection with the following invariant (recall that a \emph{parameter ideal} is an ideal generated by a system of parameters).

\begin{definition}[Parameter Degree]
For $A$ a Noetherian local ring, let its \emph{parameter degree} be the minimal  length of a homomorphic image $ A/I$, for $I$ running over all possible parameter ideals.
\end{definition}

In calculating the parameter degree, it suffices to let $I$ run over all minimal reductions of the maximal ideal of $A$. By the argument in the proof of Propositon~\ref{P:NNmult}, parameter degree and Noether Normalization degree agree for equi\ch\ complete Noetherian local rings. In fact, in view of \citep[Theorem 17.11]{Mats}, we showed the following.

\begin{corollary}
For a Noetherian local ring $A$, its multiplicity  is less than or equal to its parameter degree, which is less than or equal to its Noether Normalization degree. 

Moreover, for $A$ complete, multiplicity and parameter degree agree \iff\ $A$ is \CM.
\end{corollary}

It is not clear whether parameter degree and Noether Normalization degree coincide for a complete domain of mixed \ch. Using the homological degree introduced by Vasconcelos in \citep{VasCD,Vas,VasHD}, we obtain the following estimate.

\begin{corollary}\label{C:hommult}
For a \nd\ complete Noetherian local ring $A$, its   Noether Normalization degree is at most its homological multiplicity $e_{\text{h}}$. 

In particular, for $A$ a \nd\ Noetherian local ring, we have estimates $\textsf{cm}_1(A)\leq e_{\text{h}}$ and $\textsf{cm}^\tau_2(A)\leq (\tau+1)e_{\text{h}}$.
\end{corollary}
\begin{proof}
Recall from \citep[\S 3]{VasCD} or \citep[\S 9.5]{Vas}, the notion of \emph{homological degree} of a module. In particular, in \citep[Definition 3.23]{VasCD}, the homological multiplicity $e_{\text{h}}$ of $A$ is defined to be the homological degree of $A$ viewed as an $A$-module.  We can always find a regular subring  $S\subset A$ so that the extension is finite and such that a regular system of parameters on $S$ is part of a minimal system of generators of the maximal ideal of $A$. Such an extension can be used to calculate $e_{\text{h}}$ by \citep[Remark 3.11]{VasCD}. It follows that $e_{\text{h}}$ is also the homological degree of $A$ viewed as an $S$-module. By \citep[Proposition 4.1]{VasCD}, we then get that $e_{\text{h}}$ is a bound on the number of generators of $A$ as an $S$-module, whence a priori, on the Noether Normalization degree of $A$.

To prove the last statement, observe that $A$ and its completion have the same homological multiplicity by \citep[Theorem 3.22]{VasCD}, so that the estimates follow from Corollary~\ref{C:NNdeg}.
\end{proof}

This raises the question whether the above estimates remain true if we remove the condition that $A$ is \nd.

\begin{remark}
Perhaps it is not correct to speak of the Noether Normalization \emph{degree}, as it is not a  \emph{cohomological degree} in the sense of Vasconcelos (see \citep[Definition 3.1]{VasCD} or \citep[\S 9.5]{Vas}). For one, it is only defined for rings, not for finitely generated modules over these rings, and then only for those admitting a sort of Noether Normalization (at best, we could define it on the class of \nd\ Noetherian local rings, by setting the Noether Normalization degree of such a ring equal to the Noether Normalization degree of its completion). Nonetheless, observe the following similarities with a cohomological degree. The Noether Normalization degree coincides with classical degree (that is to say, with multiplicity) in the case the local ring $(A,\maxim)$ is \CM\ and equi\ch\ (and in certain mixed \ch\ cases). Furthermore, the Noether Normalization degree of a hyperplane section $A/hA$ is the same as the Noether Normalization degree $N$ of $A$, provided $h$ is chosen to be a regular parameter of $S$, where $S$ is a regular subring of $A$ with $\mu_S(A)=N$ (indeed, $S/hS$ is again regular and embeds in $A/hA$). Hyperplane sections $h$ as above are abundant: we can always find one such $h$ avoiding a given finite set of non-maximal prime ideals of $A$. Unfortunately, the Noether Normalization degree does not behave well with respect to elements of finite length and we have only the inequality that $N$ is at most the Noether Normalization degree of $A/\Gamma_\maxim(A)$ plus the length of $\Gamma_\maxim(A)$, where $\Gamma_\maxim(A)=\op H_\maxim^0(A)$ is the zero-th local cohomology of $A$, that is to say, the stable value of the ascending chain of ideals $\op{Ann}(\maxim^n)$. Note that for a true cohomological degree one wants equality here.

Parameter degree behaves even less well, but it can easily be defined for finitely generated modules as well. Is there therefore a formalism, similar to the 'big DEG's of Vasconcelos, that does include these degrees?
\end{remark}

\section{The Global Case}

To make the reduction to the local case, we use the Forster-Swan Theorem (see for instance \citep[Theorem 5.7]{Mats}).

\begin{citedtheorem}[Forster-Swan Theorem]
Let $A$ be a Noetherian ring and $M$ a finitely generated $A$-module. For each prime ideal $\pr$ of $A$, define
	\begin{equation*}
	f(\pr,M):= \mu_{A_\pr}(M_\pr)+\op{dim}(A/\pr).
	\end{equation*}
If $F$ is the maximum of all $f(\pr,M)$ for $\pr$ running over all prime ideals in the support of $M$, then $M$ can be generated by at most $F$ elements.
\end{citedtheorem}

\begin{corollary}\label{C:loc}
Let $A$ be a $d$-dimensional Noetherian ring and $\id$ an ideal of $A$. Let $F$ be a bound on the number of generators of each $\id A_\maxim$, where $\maxim$ runs over all maximal ideals of $A$. Then $\id$ can be generated by at most $\op{max}\{d+1,F+\op{dim}A/\id\}$ elements.
\end{corollary}
\begin{proof}
Let $\pr$ be an arbitrary prime ideal  of $A$. If $\id$ is not contained in $\pr$, then $\id A_\pr=A_\pr$ is generated by a single element, so that $f(\pr,\id)=1+\op{dim}A/\pr\leq d+1$. If $\id\subset\pr$, then $\op{dim}A/\pr\leq\op{dim}A/\id$. Choose a maximal ideal $\maxim$ of $A$, containing $\pr$. Since $\id A_\pr $ is a localization of $\id A_\maxim$, it is generated by at most $F$ elements. The assertion now follows from the Forster-Swan Theorem.
\end{proof}

We also need to study the behavior  of Noether Normalization degrees under localization and completion.

\begin{proposition}\label{P:comp}
Let $A$ be a Noetherian ring with Noether Normalization degree $N$. For every prime ideal $\pr$ of $A$, the Noether Normalization degree of the completion $\complet {A_\pr}$ of $A_\pr$ is at most $N$.
\end{proposition}
\begin{proof}
Let $B$ be a regular subring of $A$ such that $\mu_B(A)=N$ and let $\mathfrak q=\pr\cap B$.  By base change, the fiber ring $A_{\mathfrak q}/\mathfrak qA_{\mathfrak q}$ has dimension at most $N$ over the residue field $k(\mathfrak q)$. Since $\complet {A_\pr}$ is a direct summand of the $\mathfrak q$-adic completion $\complet{A_{\mathfrak q}}$ of $A_{\mathfrak q}$ by \citep[Theorem 8.15]{Mats}, we get that
	\begin{equation*}
	\op{dim}_{k(\mathfrak q)}(\complet {A_\pr}/\mathfrak q\complet {A_\pr})\leq \op{dim}_{k(\mathfrak q)}(\complet {A_{\mathfrak q}}/\mathfrak q\complet {A_{\mathfrak q}})=\op{dim}_{k(\mathfrak q)}( A_{\mathfrak q}/\mathfrak qA_{\mathfrak q})\leq N.
	\end{equation*}
In particular, $\complet {A_\pr}$ is generated as a $\complet{B_{\mathfrak q}}$-module by at most $N$ elements, by \citep[Theorem 8.4]{Mats}. Since $B_{\mathfrak q}$ is regular, whence also its completion, we showed that $\complet{A_\pr}$ has Noether Normalization degree at most $N$.
\end{proof}

\begin{theorem}\label{T:eqch}
Let $A$ be a $d$-dimensional Noetherian ring of finite   Noether Normalization degree $N$. We have estimates   $\textsf{cm}_1(A)\leq N+d-1$ and $\textsf{cm}_2^\tau(A)\leq (\tau+1)N+d-2$ (except when $N=\tau=1$, in which case possibly $d+1$ generators are needed).
\end{theorem}
\begin{proof}
Let $\id$ be a \CM\ ideal of $A$ of height at most $2$. By Corollary~\ref{C:loc}, if we find a bound $F$ on the number of generators of $\id A_\maxim$ in each localization with respect to a maximal ideal $\maxim$, then $\id$ itself can be generated by $F+\op{dim}A/\id$ elements. The statement therefore follows from Corollary~\ref{C:NNdeg} and Proposition~\ref{P:comp}, since we may take $F=N$. One just needs to observe that in the height one case $d+1\leq d+N-1$, except when $N=1$. However, $N=1$ means that $A$ is regular, whence in particular a \UFD, so that any height one ideal is principal. In the height two case, the given bound is at least $d$ and only in the indicated case  $N=\tau=1$ it is equal to it.\footnote{In the latter case,  $A$ is regular and $\id$ is a Gorenstein ideal, whence locally a complete intersection. If $A$ is a polynomial ring over some subring, then the EE-Conjecture proven by \citet{MK}, states that we may drop the contribution of $f(\pr,\id)$ for all minimal primes $\pr$ of $A$ in the bound in the  Forster-Swan Theorem, yielding therefore in this case the correct bound $d$.  Is it true that every height two locally complete intersection in a $d$-dimensional regular ring is generated by at most $d$ elements?}
\end{proof}

\subsubsection*{Proof of Theorem~\ref{T:A}}
If $S$ and $A$ denote the affine algebras of $X$ and $Y$ respectively, then our assumptions imply that $S\subset A$ is finite with $S$ regular. By definition, the  degree $f$ of $Y\to X$ is the maximal number of points in a closed fiber. In other words, $f$ is the maximum of the dimensions 
	\begin{equation*}
	f(\maxim):=\op{dim}_{S/\maxim} (A_\maxim/\maxim A_\maxim),
	\end{equation*}
where $\maxim$ runs over all maximal ideals of $S$. Since $f(\maxim)$ is equal to the minimal number of generators of $A_\maxim$ over $S_\maxim$, a similar argument as in the proof of Corollary~\ref{C:loc} shows that $A$ is generated as an $S$-module by at most $f+d$ elements. The stated bound now follows by letting $N:=f+d$ in Theorem~\ref{T:eqch}.

If $Y$, or equivalently $A$, is \CM, then the multiplicity of $A_{\mathfrak n}$, for $\mathfrak n$ a maximal ideal of $A$, is equal to its Noether Normalization degree, which is at most $f(\maxim)\leq f$, for $\maxim:=\mathfrak n\cap S$ (apply Propositions~\ref{P:NNmult} and \ref{P:comp} to the completion of $A_{\mathfrak n}$). By the Main Theorem, we may take $F:=f$ in Corollary~\ref{C:loc}, so that we get the bound $2f+d-2$ on the number of generators of a height two Gorenstein ideal in $A$.
\qed

\begin{example}
If $X$ is an affine subspace over some other scheme, then by the EE-Conjecture (see the previous footnote), me may take $N:=f+d-1$, so that in this case we only require $3d+2f-4$ hypersurfaces. For instance, if $Y$ is an $f$-fold covering of $\mathbb A_K^3$, with $K$ a field, then every local Gorenstein curve in $Y$ is the (ideal-theoretic) intersection of at most $2f+5$ surfaces (at most $2f+1$, if $Y$ is \CM). 
\end{example}

\section{The Affine Case}

Affine rings, that is to say, finitely generated algebras over a field, have the property that their Noether Normalization degree is finite. In fact, we have the following sharper result.

\begin{theorem}\label{T:aff}
For each $n$, there exists a bound $\bound En$ with the following property. If $A$ is an affine ring of the form $\pol KX/\rij fs\pol KX$ with $K$ a field, $X$ a set of at most $n$ variables  and $f_i$ polynomials of degree at most $n$, then the Noether Normalization degree of $A$ is at most $\bound En$.

In particular,  we have estimates  $\textsf{cm}_1(A)\leq \bound En+n$ and  $\textsf{cm}_2^\tau(A)\leq (\tau+1)\bound En+n$.
\end{theorem}
\begin{proof}
The second statement follows from the first and Theorem~\ref{T:eqch}, as $A$ has dimension at most $n$. 

To prove the first statement, one just needs to observe that Noether Normalization can be carried out algorithmically from the $f_i$. The key idea is to make a change of variables so that  one of the $f_i$ becomes monic in some variable. If $K$ is infinite, this can be done by a linear change of variables; in the general case, we can still control the degree of this new equation (see \citep[\S A.5]{Vas} for details). Assume therefore that all $f_i$ have degree at most $n'$ and that $f_1$ is monic in $X_1$ of degree $n'$, where $n'$ only depends on $n$. Hence $\pol KX/f_1\pol KX$ is generated by $1,X_1,\dots,X_1^{n'-1}$ over $\pol K{X_2,\dots,X_n}$. Let 
	\begin{equation*}
	I_1:=\rij fs\pol KX\cap \pol K{X_2,\dots,X_n}
	\end{equation*}
 and put $A_1:=\pol K{X_2,\dots,X_n}/I_1$. It follows that $A_1\subset A$ is a finite extension, generated by at most $n'$ elements. By  \citep[Theorem 2.6]{SchBC}, the ideal $I_1$ is generated by polynomials of degree at most $n''$, where $n''$ depends only on $n'$, whence only on $n$. Therefore, by an inductive argument, $A_1$ admits a Noether Normalization $\pol KY\subset A_1$ generated by at most $n'''$ elements, where $n'''$ depends only on $n''$, whence only on $n$. From   the composition $\pol KY\subset A_1\subset A$ we see that $A$ is generated as a $\pol KY$-module by at most $n'''n''$ elements, a number only depending on $n$.
\end{proof}

In fact, the above bounds hold uniformly in families.

\begin{corollary}\label{C:fam}
Let $s\colon W\to V$ be map of finite type between schemes of finite type over some field. There exists a bound $\bound{CM}s$, such that for each $x\in V$ and each codimension two \CM\ subscheme $F$ of $\inverse sx$ of type $\tau$, the ideal of $F$ is generated by at most $(\tau+1)\bound{CM}s$ elements. If $F$ has codimension one, then at most $\bound{CM}s$ generators suffice.
\end{corollary} 
\begin{proof}
Taking a finite affine covering, we may reduce to the case that $V=\op{Spec}A$ and $W=\op{Spec}B$ are affine, so that $s$ corresponds to a $K$-algebra \homo\ $A\to B$. We can write $A$ as a homomorphic image of a polynomial ring $\pol KX$ modulo an ideal $\rij fs\pol KX$. Moreover, we can find $g_i\in\pol K{X,Y}$, such that $B\iso \pol AY/\rij gs\pol AY$. Let $n$ be the maximum of the number of variables and of the degrees of the $f_i$ and the $g_i$. If $\pr$ denotes the prime ideal of $A$ corresponding to the point $x\in V$, then the coordinate ring of the fiber $\inverse sx$ is $B\tensor_A k(\pr)$, where $k(\pr):=A_\pr/\pr A_\pr$ is the residue field of $\pr$. It follows that $B\tensor_Ak(\pr)\iso\pol{k(\pr)}X/\rij gs  \pol{k(\pr)}X$. By Theorem~\ref{T:aff}, the ideal of $B\tensor_Ak(\pr)$ defining $F$ is generated by at most $(\tau+1)\bound En+n$ elements, where $\bound En$ is as in that Theorem. The height one case follows by a similar argument.
\end{proof}

The proof actually shows that $\bound{CM}s$ only depends on the degrees of the polynomials involved in defining $V$, $W$ and $s$.

\begin{corollary}\label{C:quot}
For each $n$, there exists a bound $\bound Nn$ with the following property. Let $A$ be an affine ring of the form $\pol KX/\rij fs\pol KX$ with $K$ a field, $X$ a set of at most $n$ variables  and $f_i$ polynomials of degree at most $n$. Let $\id$ be a height $h$ \CM\ ideal of $A$ for which the type of  $A/\id$ is $\tau$. If $\id$ contains an height $h-2$ ideal $I:=\rij gsA$, with the $g_i$ of degree at most $n$, then $\id$ can be generated by at most $(\tau+1)\bound Nn$ elements. If $I$ has height $h-1$, then  at most $\bound Nn$ generators suffice.

In particular, every height three ideal $\id$ of $A$ which contains the image of a polynomial of degree at most $n$ not belonging to any minimal prime of $A$ and for which $A/\id$ is Gorenstein, can be generated by at most $2\bound Nn$ elements.
\end{corollary}
\begin{proof}
The second statement is a special case of the first, with $h=3$ and $\tau=1$. Let $A$, $\id$ and $I$ be as in the first statement. It follows that there is a bound $\bound{N'}n$ on the number of generators of $I$, only depending on $n$.  Let $B:=A/I$, so that $B$ is also a homomorphic image of a polynomial ring in at most $n$ variables by an ideal generated by polynomials of degree at most $n$. Applying Theorem~\ref{T:aff} to the height two \CM\ ideal $\id B$, it follows that $\id B$ is generated by at most $(\tau+1)\bound En$ elements. Therefore, $\id$ is generated by at most $(\tau+1)\bound En+\bound{N'}n$ elements. 
\end{proof}


\begin{thebibliography}{14}
\expandafter\ifx\csname natexlab\endcsname\relax\def\natexlab#1{#1}\fi

\bibitem[Boraty\'nski et~al.(1979)Boraty\'nski, Eisenbud, and Rees]{BER}
M.~Boraty\'nski, D.~Eisenbud, and D.~Rees.
\newblock On the number of generators of ideals in local {C}ohen-{M}acaulay
  rings.
\newblock {\em J. Algebra}, 57:\penalty0 77--81, 1979.

\bibitem[Bourbaki(1983)]{Bour83}
N.~Bourbaki.
\newblock {\em Alg{\` e}bre Commutative}, volume 8-9.
\newblock Masson, New York, 1983.

\bibitem[Bruns and Herzog(1993)]{BH}
W.~Bruns and J.~Herzog.
\newblock {\em Cohen-{M}acaulay Rings}.
\newblock Cambridge University Press, Cambridge, 1993.

\bibitem[Matsumura(1986)]{Mats}
H.~Matsumura.
\newblock {\em Commutative Ring Theory}.
\newblock Cambridge University Press, Cambridge, 1986.

\bibitem[Mohan-Kumar(1978)]{MK}
N.~Mohan-Kumar.
\newblock On two conjectures about polynomial rings.
\newblock {\em Invent. Math.}, 46:\penalty0 225--236, 1978.

\bibitem[Sally(1978)]{Sally}
J.~Sally.
\newblock {\em Numbers of Generators of Ideals in Local Rings}.
\newblock Lect. Notes in Pure and Applied Mathematics. Marcel Dekker, Inc.,
  1978.

\bibitem[Sally and Vasconcelos(1974)]{SV}
J.~Sally and W.~Vasconcelos.
\newblock Stable rings.
\newblock {\em J. Pure Appl. Algebra}, 4:\penalty0 319--336, 1974.

\bibitem[Schoutens(2000)]{SchBC}
H.~Schoutens.
\newblock Bounds in cohomology.
\newblock {\em Israel J. Math.}, 116:\penalty0 125--169, 2000.

\bibitem[Schoutens(2003)]{SchCMId}
H.~Schoutens.
\newblock Number of generators of a {C}ohen-{M}acaulay ideal.
\newblock {\em J. Algebra}, 259:\penalty0 235--242, 2003.

\bibitem[Valla(1981)]{VallaGen}
G.~Valla.
\newblock Generators of ideals and multiplicities.
\newblock {\em Comm. Algebra}, 9:\penalty0 1541--1549, 1981.

\bibitem[Vasconcelos(1998{\natexlab{a}})]{VasCD}
W.~Vasconcelos.
\newblock Cohomological degrees of graded modules.
\newblock In J.~Elias, J.M. Giral, R.M. Mir\'o-Roig, and S.~Zarzuela, editors,
  {\em Six lectures on Commutative Algebra}, volume 166 of {\em Progress in
  Math.}, pages 345--390. Birkh\"auser, 1998{\natexlab{a}}.

\bibitem[Vasconcelos(1998{\natexlab{b}})]{Vas}
W.~Vasconcelos.
\newblock {\em Computational Methods in Commutative Algebra and Algebraic
  Geometry}.
\newblock Springer-Verlag, Berlin, 1998{\natexlab{b}}.

\bibitem[Vasconcelos(1998{\natexlab{c}})]{VasHD}
W.~Vasconcelos.
\newblock The homological degree of a module.
\newblock {\em Trans. Amer. Math. Soc.}, 350:\penalty0 1167--1179,
  1998{\natexlab{c}}.

\bibitem[Vasconcelos(2002)]{VasMultGen}
W.~Vasconcelos.
\newblock Multiplicities and the number of generators of {C}ohen-{M}acaulay
  ideals.
\newblock manuscript, 2002.

\end{thebibliography}
\end{document}